\documentclass[12pt]{article}
\newtheorem{prelemmaa}{{\bf LEMMA}}

\newtheorem{prelem}{{\bf THEOREM}}

\newtheorem{preque}{{\bf QUESTION}}

\newtheorem{theorem}{THEOREM}

\newtheorem{prelemma}{LEMMA}

\newenvironment{lemma}{\begin{prelemma}{\hspace{-0.5
               em}}}{\end{prelemma}}
\newtheorem{preproof}{{\bf PROOF.}}

\newenvironment{proof}[1]{\begin{preproof}{\rm
               #1}\hfill{\rule[-0.5mm]{2mm}{2mm}}}{\end{preproof}}

\newtheorem{preproposition}{{PROPOSITION}}

\newtheorem{preremark}{REMARK}

\newtheorem{precorollary}{{COROLLARY}}

\newtheorem{precorn}{{COROLLARY}}

\newtheorem{predefinition}{DEFINITION}

\newtheorem{preexample}{EXAMPLE}

\newtheorem{preconjecture}{{CONJECTURE}}

\newenvironment{conjecture}{\begin{preconjecture}{\hspace{-0.5
                em}}}{\end{preconjecture}}
\newtheorem{pretheo}{{\bf THEOREM}}

\newenvironment{theo}{\begin{pretheo}{\hspace{-0.5
               em}{\bf}}}{\end{pretheo}}
\def\newpic#1{}
\date{}
\begin{document}
\title{\bf On the non-existence of some Steiner $t$-$(v,k)$ trades of
certain volumes}
\author{
{\sc Mehri Asgari AND Nasrin Soltankhah} \footnote{Corresponding
author: E-mail: soltan@alzahra.ac.ir.}
 \\[5mm]
Department of Mathematics\\ Alzahra University  \\
Vanak Square 19834 \ Tehran, I.R. Iran }
%
\maketitle
\begin{abstract}
Mahmoodian and Soltankhah $\cite{MMS}$ conjectured that there does
not exist any $t$-$(v,k)$ trade of volume $s_{i}< s <s_{i+1}$,
where $s_{i}=2^{t+1}-2^{t-i},\  i=0,1,\ldots ,t-1$. Also they
showed that the conjecture is true for $i=0$. In this paper we
prove the correctness of this conjecture for Steiner trades.
 \end{abstract}
%
\hspace*{-2.7mm} {\bf KEYWORDS:} {  \sf trade, Steiner trade,
volume }
%
\section{Introduction} 
\setcounter{preproposition}{0} \setcounter{precorollary}{0}
\setcounter{prelemma}{0} \setcounter{preexample}{0} Let $0<t<k<v $
be three natural numbers and let $X$ be a $v$-set. For every $i,\
0\leq i\leq v $, the set of all $i$-subsets of $X$  will be
denoted by $P_i(X)$, and also for any $i>1 $, we will denote the
set $\{x_1,x_2,\ldots,x_i\} $ by $x_1x_2\ldots x_i$. The elements
of $P_k(X)$ are called \textbf{blocks}. A $t$-$(v,k)$
\textbf{trade}
 $T=\{T_1,T_2\}$ consists of two disjoint
collections of blocks, $T_1$ and $T_2$, such that every element
of $P_t(X)$ is contained  in the same number of blocks in $T_1$
and $T_2$. For simplicity, the notation of $t$-trade is commonly
used in this manuscript. Let $T=\{T_1,T_2\}$ be a $t$-trade.
Clearly, $|T_1|=|T_2|$ and $|T_1|$ is called the \textbf{volume}
of the $T$ and is denoted by \rm {vol(T)} and sometimes we denote
vol(T) by \rm {s}, if it causes no confusion. The subset of $X$
which is covered by $T_1$ and $T_2$ is called the
\textbf{foundation} of T and is denoted by \rm {found(T)}.
Repeated blocks in $T_1(T_2)$ are allowed.

Trade with no repeated block is simple. A $t$-trade is called
\textbf{Steiner} $t$-$(v,k)$ \textbf{trade} if every element of
$P_t(X)$ appears in at most one block of $T_1(T_2)$. It has been
shown in $\cite{{ graham},{HL}}$ that in every $t$-$(v,k)$ trade,
$|$found$(T)|\geq k+t+1$, vol$(T)\geq2^t$. A $t$-trade $T$ with
vol$(T)=2^t$ and $|$found$(T)|= k+t+1$ is called
\textbf{minimal}. A minimal $t$-$(v,k)$ trade is unique, up to
isomorphism, and can be cast in the following form
$T=(x_1-x_2)(x_3-x_4)\ldots(x_{2t+1} - x_{2t+2})x_{2t+3} \ldots
x_{k+t+1}$, where $x_i \in $ found$(T)$. After a formal
multiplication, the terms with plus (minus) signs are to be
considered as block of $T_1(T_2)$. Let $T=\{T_1,T_2\}$ be a
$t$-$(v,k)$ trade of vol$(T)=s$ and $x,y\in$ found$(T)$. Then the
number of blocks in $T_1(T_2)$ which contains $x$ is denoted by
$r_x$, and the number of blocks containing $\{x,y\}$ (for
$t\geq2)$  is denoted by $\lambda_{xy}$. The set of blocks in
$T_1(T_2)$ which contains $x \in$ found$(T)$ is denoted by
$T_{1x}(T_{2x})$ and the set of remaining blocks by
$T^\prime_{1x}(T^\prime_{2x})$. It has been shown $\cite{HL}$ that
if $r_x < s$, then $T_x=\{T_{1x},T_{2x}\}$, is a $(t-1)$-$(v,k)$
trade of vol$(T_x)=r_x$, and furthermore,
$T^\prime_x=\{T^\prime_{1x},T^\prime_{2x}\}$ is a
$(t-1)$-$(v-1,k)$ trade of vol$(T^\prime_x)=s-r_x$. If we remove
$x$ from the blocks of $T_x$ then the result will be a
$(t-1)$-$(v-1,k-1)$ trade which is called a derived trade of $T$
and is denoted by $D_x T=\{(D_x T)_1, (D_x T)_2\}$. It is easy to
show that if T is a Steiner trade, then its derived trade is also
a Steiner trade. If $T=\{T_1,T_2\}$ and $T^*=\{T^*_1,T^*_2\}$ are
two $t$-$(v,k)$ trades, then we define $T+T^*=\{T_1\bigcup
T^*_1,T_2\bigcup T^*_2 \}$ and $T-T^*=\{T_1\bigcup
T^*_2,T^*_1\bigcup T_2\}$. Note that the blocks that appear in
both sides are omitted. It is easy to see that $T \pm T^*$ are
also $t$-$(v,k)$ trades. Let $T$  be a $t$-$(v,k)$ trade and $T
\neq T_x +T_y $ (for $x,y \in$ found$(T))$ then $T-(T_x+T_y)$ is a
$(t-1)$-$(v,k)$ trade with volume $s-(r_x+r_y)+2\lambda_{xy}$
$\cite{HL}$. Hwang $\cite{HL}$ has shown that there is no
$2$-$(v,k)$ trade of volume $s=5$, and as a generalization, she
has also shown that $t$-trades with $s=2^t+1$ do not exist,
Mahmoodian and Soltankhah $\cite{MMS}$ have shown that:
 \begin{theo}~{\rm \cite{MMS}}. There does not exist any $t$-$(v,k)$ trade of volume $s$ such that $2^t<s<2^t + 2^{t-1}$.
\end{theo}
 Also in \cite{MMS}, the
following conjecture has been stated:
 \begin{conjecture}~{\rm \cite{MMS}}.
  There does not exist any $t$-$(v,k)$ trade of {\rm volume} $s$ such that
$2^{t+1} - 2^{t-i}<s<2^{t+1} - 2^{t-i-1},\ \ i=0,1,\ldots,t-1$.
\end{conjecture}
In fact Theorem A is the answer for the conjecture for $i=0$.
Gray and Ramsay $\cite{Gray}$, showed that no $3$-$(v,4)$ trade of
$s=13$ exists, and they also proved a generalization of this: For
$t\geq 3$, $t$-$(v,t+1)$ trades of volume $s=2^t+ 2^{t-1}+1$  do
not exist. Of course this is a partial solution for $i=1$ of the
conjecture. Prior to this proof, in $\cite{hart}$ it is shown
that there does not exist any Steiner $3$-$(v,4)$ trade of $s=13$.
Hoorfar and G.B. Khosrovshahi $\cite{khosro}$ established the
correctness of this conjecture for Steiner trades, for $i=1$.
 In this paper, we prove the correctness of the conjecture for Steiner
 trades,   for $i=0,\cdots,t-1$. In other words we show that the
 conjecture is true for Steiner trades .
\section{Some necessary results } 
\setcounter{theorem}{0} \setcounter{preproposition}{0}
\setcounter{precorollary}{0} \setcounter{prelemma}{0}
\setcounter{preexample}{0}

The  following results are  useful in our discussion.
\begin{lemma}~{\rm \cite{khosro}}.
\label{l1} In every Steiner $t$-trade $T=\{T_1,T_2\}$ with  $k =
t+1$, volume $s$ and for every $x\in $found$(T)$, we have $r_x\leq
\frac{1}{2} s$.

\end{lemma}
\begin{theorem}~{\rm \cite{khosro}}.\label{c2}
Let $T=\{T_1,T_2\}$ be a Steiner $t$-$(v,k)$ trade with $k> t+1$
and of volume $s$, then $ s \geq (t-1)2^t +2$.

\end{theorem}
%

\section{Main result}  
In this section we investigate the main theorem in two cases:
$k=t+1$ and $k>t+1$. Note that as a result of Theorem~\ref{c2} ,
it follows that for $t=2$, we have $s\geq6$ and for $t>2$, we
have $s\geq 2^{t+1}+2$. Therefore Theorem~\ref{c2} establishes
the truth of the Conjecture for any Steiner trade with $k> t+1$.
Therefore the case $k=t+1$ remains to prove, which we discuss it
in the following Theorem.

\begin{theorem}.
 For $t\geq 3$, there does not exist any Steiner
$t$-$(v,t+1)$ trade of volume $s$ such that
$2^{t+1}-2^{t-i}<s<2^{t+1}-2^{t-i-1},\  i=0,1,\ldots t-1 $.

\end{theorem}
\begin{proof}
{Induction on $i$. For $i=0$, this is true by Theorem A. Suppose
that the theorem is correct for values smaller than
$i\ (i>1)$ and we have to establish it for $i$.\\
 Now induction on $t$.\\
 The minimum possible value on t is equal to $i+2$.
First we show, the nonexistence of these $t$-trades.\\
 If $t=i+2$, then $s=2^{t+1}-3$.
Suppose that the statement is not correct and there exists a
Steiner $t$-$(v,t+1)$ trade $T$ of volume $s=2^{t+1}-3$ and
$|$found$(T)|=f$. Then we derive a contradiction.\\
 Consider $x\in$ found$(T)$, since $D_x T $ is a Steiner $(t-1)$-trade
with $k=t$ and volume $s=r_x$, therefore $r_x\geq 2^{t-1} $. On
the other hand by Lemma 2.1, we
have $r_x\leq \frac{s}{2}$.\\
$ \Rightarrow r_x\leq \frac{s}{2}=\frac{2^{t+1}-3}{2}\leq 2^t -
2$.\\
 Therefore $2^{t-1}\leq r_x \leq 2^t-2$.  But by induction assumption,
only values remain to be checked are: $r_x =2^t -2^{t-u}\
(u=1,\ldots ,i+1)$.
\vspace{.2cm}\\
 \textbf{{Case A.1:}}  there exists $x\in$ found$(T)$ such
that $r_x=2^{t-1}$ then,\\
 \textbf{\textrm(a)} suppose that there
exists $y\in$ found$(T)$ such that $\lambda_{xy}=0 $. this
\textsl{leads} to a contradiction. To see this we look at $T-(T_x
+T_y)$ which is a $(t-1)$-$(v,k')$ trade of volume $
s^\prime=s-(r_x+r_y)
=2^{t+1}-3-2^{t-1}-r_y.$\\
$\Rightarrow 2^t+2^{t-1}-3-2^t+2\leq s^\prime
\leq 2^t+2^{t-1}-3-2^{t-1}$;\\
$\Rightarrow 2^{t-1}-1\leq s^\prime \leq 2^t-3$.\\
 If $k' >t $, by Theorem 2.1 is
impossible.
 If $k'=t$ since $s^\prime $ is odd and by induction assumption,
 $s\prime $ can not take a possible value.\\
 \textbf{\textrm{(b)}} suppose that for every $y\in$ found$(T)$, $y\neq
x$
 , we have $\lambda_{xy}\neq 0$ then $|$found$(T_x)|=|$found$(T)|=f$.
  Now, since $T_x$ is a $(t-1)$-trade with $k=t+1$ and volume
$r_x=2^{t-1}$ such that $x$ appears in all its blocks. Must be of
the form
\begin{center}
$T_x=(y_1-y_2)\ldots (y_{2t-1}-y_{2t})x.$
\end{center}
 We note that $f=2t+1$. This is a contradiction since in
$T$ we have
\begin{center}
$ f \geq  k+t+1=2t+2$.
\end{center}
\textbf{{Case A.2:}} for every $x\in$ found$(T):\  r_x\geq
2^t-2^{t-2}.
$\\
\textbf{\textrm{(a)}} again suppose that there exists $y\in$
found$(T);$ \ $y\neq x $ such that  $\lambda_{xy}=0. $ Then the
volume $s^\prime$ of the $(t-1)$-trade $T-(T_x+T_y)$ is:
$s^\prime=s-(r_x+r_y)=2^{t+1}-3-(r_x+r_y)$. With considering the
values $r_x$ and $r_y$ we have
\begin{center}
$1\leq s^\prime\leq 2^{t-1}-3$ \\
\end{center}
since minimum volume from $(t-1)$-trade is $2^{t-1}$ it is a contradiction.\\
 \textbf{\textrm{(b)}} now for every $x,y \in$ found$(T)$ we have
 $\lambda_{xy}
 \neq0 .$\\
 Therefore:\\
 if $ r_x=2^t-2^{t-u}, \hspace{2 cm}u=2,\ldots ,i+1 $;\\
       $\lambda_{xy}=2^{t-1}-2^{t-l},\hspace{2 cm}  l=2,\ldots, u+1.
$\\
\textbf{i:} $\exists y\in$ found$(T),\ \lambda_{xy} =
2^{t-1}-2^{t-l},\hspace{1 cm} l=2,\ldots, u$.\\
In this case volume $(t-1)$-trade $T-(T_x+T_y)$ is:
$s^\prime=s-(r_x+r_y)+2\lambda_{xy}$
\begin{center}
$\Rightarrow 2^{t-1}+1\leq s^\prime \leq 2^t-3$.
\end{center}
 Since $s^\prime $ is odd and by induction assumption,
 $s^\prime$ can not take a possible value.\\
\textbf{ii:} $\forall y\in$ found$(T),\
\lambda_{xy}=2^{t-1}-2^{t-u-1}$.\\
 In this case to reach a contradiction, counting
the pairs $(y,B); y\in B\in (D_xT) $ in two ways we obtain
$$
 \sum_{y\in {\rm found}(D_xT)}\lambda_{xy}=t.{\rm vol}(D_xT).$$
Since $|$found$(D_xT)|=f-1$, therefore
\begin{center}
$(f-1)(2^{t-1}-2^{t-u-1})=t(2^{t}+2^{t-u})$ \\
$\Rightarrow f=2t+1,$
\end{center}
and this is in contradiction with $ f\geq k+t+1=2t+2$.\\

Now suppose by induction assumption that the theorem is correct
for values smaller than $t(t>i+2)$ and we show it for $t$. Suppose
that the statement is not correct and there exists a Steiner
$t$-$(v,t+1)$ trade T of volume $s=2^t+2^{t-1}+\cdots +2^{t-i}+j$
where $ 0<j<2^{t-1-i}.
$ Then we derive a contradiction.\\
Consider $x\in$ found$(T)$. Since $D_xT$ is a Steiner
$(t-1)$-trade with $k=t$ and volume $s=r_x$, therefore $r_x\geq
2^{t-1}$ and on  the other hand we have $r_x\leq \frac{s}{2}$.
Now by induction
assumption:\\
$r_x=2^t-2^{t-u}$ where $u=1,2,\ldots,i+1. $
\vspace{.2cm}\\
 \textbf{{Case B.1:} } there exists $x\in$ found$(T)$ such
that $r_x=2^{t-1}$ then\\
 \textbf{\textrm(a)} suppose that there
exists $y\in$ found$(T)$ such that $\lambda_{xy}=0$. For
$(t-1)$-trade $T-(T_x+T_y)$
\begin{center}
$s^\prime
=s-(r_x+r_y)=2^{t+1}-2^{t-i}+j-(2^{t-1}+2^t-2^{t-u})=2^{t-1}+2^{t-u}-2^{t-i}+j.$
\end{center}
$\begin{array}{l} {\rm if}\hspace{1 cm}u=1 \hspace{3
cm}s^\prime=2^t-2^{t-i}+j
 \hspace{3 cm} 0<j<2^{t-i-1}\\
\end{array}$\\
$\begin{array}{l} {\rm if}\hspace{0.5 cm}2\leq u\leq i
\hspace{2cm}$
$2^{t-1}<s^\prime<2^{t-1}+2^{t-2}\\
\end{array}$\\
$\begin{array}{l}
 {\rm if}\hspace{1 cm} u=i+1 \hspace{2.5 cm}
s^\prime=2^{t-1}-2^{t-i-1}+j<2^{t-1}-2^{t-i-1}+2^{t-i-1}<2^{t-1}.
\end{array}$\\
In the first case according to induction assumption  for $(i-1)$,
there does not exist any $(t-1)$-trade of this volume.
\vspace{.2cm}\\
 In the second case  according to induction assumption
 for $i=0$ there does not exist any  $(t-1)$-trade of
this volume.
\vspace{.25cm}\\
 And ultimately the last case is in contradiction
with minimum volume from $(t-1)$-trade.
\vspace{.2 cm}\\
 \textbf{\textrm(b)} suppose that for every $y\in$ found$(T)$, $y\neq x$
 we have $\lambda_{xy}\neq 0.$\\
 The
non-existence of T can be proved the same way as in the case
$t=i+2.$
 \vspace{.2cm}\\
 \textbf{{Case B.2:}} for every $x\in$ found$(T)\  r_x\geq 2^t-2^{t-2}. $\\
 \textbf{\textrm{(a)}}
 Again suppose that there exists $y\in$ found$(T) ;\ y\neq x $
 such that $\lambda_{xy}=0. $\\
 For $(t-1)$-trade $T-(T_x+T_y)$\\
$ s^\prime
=s-(r_x+r_y)=2^{t+1}-2^{t-i}+j-(2^t-2^{t-u}+2^t-2^{t-u^\prime})
\hspace{0.2 cm}$
where $u,u^\prime \geq2$\\
$\Rightarrow s^\prime =2^{t-u}+2^{t-u^\prime}-2^{t-i}+j\leq
2^{t-2}+2^{t-2}-2^{t-i}+j \leq 2^{t-1}-2^{t-i}+j$\\
by assumption $0<j<2^{t-i-1}$ we have $s\prime <2^{t-1}$ which
this is in contradiction with minimum volume from $(t-1)$-trade.\\
 \textbf{\textrm{(b)}} now for every $x,y \in$ found$(T),\ r_x\leq r_y $
we have
 $\lambda_{xy}
 \neq0 .$\\
 if $ r_x=2^t-2^{t-u},$ \hspace{1 cm}  $  u=2,\ldots ,i+1 $\\
      $\lambda_{xy}=2^{t-1}-2^{t-l}, \hspace{0.5 cm} l=2,\ldots, u+1$
\\
       $r_y=2^{t}-2^{t-u^\prime}, \hspace{1 cm} u^\prime=u,\ldots
,i+1$.\\
 \textbf{i:} $\exists y\in$ found$(T),$ \ $\lambda_{xy} =
2^{t-1}-2^{t-l}, \ l=2,\ldots, u$.\\
  In this case for $(t-1)$-trade $T-(T_x+T_y)$ have:
\begin{center}
$s^\prime=2^t+2^{t-u}+2^{t-u^\prime}-2^{t-i}-2^{t-l+1}+j$\\
where  $u=2,\ldots ,i+1$;\\
\hspace{1cm} $l=2,\ldots ,u$;\\
\hspace{2.5cm} $u^\prime=u, \ldots, i+1$.\\
\end{center}
 Therefore $l\leq u\leq u^\prime $
 {we have four cases}:
 \vspace{.22 cm}\\
 \textbf{1)} $u=u^\prime=l$ \hspace{1 cm} $\Rightarrow s^\prime
 =2^t-2^{t-i}+j$.
\vspace{.22 cm}\\
 $(s^\prime=2^t+2^{t-u}+
2^{t-u^\prime}-2^{t-i}-2^{t-l+1}+j=2^t+2^{t-u}+
2^{t-u}-2^{t-i}-2^{t-u+1}+j=2^t-2^{t-i}+j)$
\vspace{.1 cm}\\
\textbf{2)} $u=l<u^\prime\leq i$  \hspace{1 cm} $\Rightarrow
2^t-2^{t-l}<s^\prime<2^t-2^{t-l-1}$.
\vspace{.22 cm}\\
 $(s^\prime=2^t+2^{t-u}+
2^{t-u^\prime}-2^{t-i}-2^{t-l+1}+j=2^t+2^{t-l}+
2^{t-u^\prime}-2^{t-i}-2^{t-l+1}+j=2^t-2^{t-l}+\underbrace{(2^{t-u^\prime}-2^{t-i})}_{\geq0}+\underbrace{j}_{>0}>2^t-2^{t-l})$
\vspace{.15 cm}\\
 $(s^\prime=2^t+2^{t-u}+
2^{t-u^\prime}-2^{t-i}-2^{t-l+1}+j=2^t+2^{t-l}+
2^{t-u^\prime}-2^{t-i}-2^{t-l+1}+j=2^t-2^{t-l}+2^{t-u^\prime}-2^{t-i}+j=2^t-2^{t-l-1}-2^{t-l-1}
+2^{t-u^\prime}-2^{t-i}+j=2^t-2^{t-l-1}+\underbrace{(2^{t-u^\prime}-2^{t-l-1})}_{\leq0}+\underbrace{(j-2^{t-i})}_{<0}<2^t-2^{t-l-1})$
\vspace{.22 cm}\\
 \textbf{3)} $u=l<u^\prime=i+1$  \hspace{1 cm} $\Rightarrow
2^t-2^{t-l+1}<s^\prime<2^t-2^{t-l}$.
\vspace{.22 cm}\\
 $(s^\prime=2^t+2^{t-u}+
2^{t-u^\prime}-2^{t-i}-2^{t-l+1}+j=2^t+2^{t-l}+
2^{t-i-1}-2^{t-i}-2^{t-l+1}+j=2^t-2^{t-l+1}+
\underbrace{(2^{t-l}-2^{t-i-1})}_{>0}+\underbrace{j}_{>0}>2^t-2^{t-l+1})$
\vspace{.15 cm}\\
$(s^\prime=2^t+2^{t-u}+
2^{t-u^\prime}-2^{t-i}-2^{t-l+1}+j=2^t+2^{t-l}+
2^{t-i-1}-2^{t-i}-2^{t-l+1}+j=2^t-2^{t-l}+2^{t-i-1}+j=2^t-2^{t-l}+\underbrace{(j-2^{t-i-1})}_{<0}
<2^t-2^{t-l})$
\vspace{.22 cm}\\
 \textbf{4)} $l<u\leq
u^\prime$\hspace{1 cm} $\Rightarrow
2^t-2^{t-l+1}<s^\prime<2^t-2^{t-l}$.
\vspace{.22 cm}\\
$(s^\prime=2^t+2^{t-u}+
2^{t-u^\prime}-2^{t-i}-2^{t-l+1}+j\geq2^t+2^{t-i-1}+
2^{t-i-1}-2^{t-i}-2^{t-l+1}+j=2^t-2^{t-l+1}+\underbrace{j}_{>0}>2^t-2^{t-l+1})$
\vspace{.15 cm}\\
 $(s^\prime=2^t+2^{t-u}+
2^{t-u^\prime}-2^{t-i}-2^{t-l+1}+j=2^t-2^{t-l}-
2^{t-l}+2^{t-u}+2^{t-u^\prime}-2^{t-i}+j\leq2^t-2^{t-l}-
2^{t-l}+2^{t-u}+2^{t-u}-2^{t-i}+j=2^t-2^{t-l}-2^{t-l}+2^{t-u+1}-2^{t-i}+j=
2^t-2^{t-l}+\underbrace{(2^{t-u+1}-2^{t-l})}_{\leq0}+\underbrace{(j-2^{t-i})}_{<0}<2^t-2^{t-l})$\\

All of the above cases are in contradiction with induction
assumption.\\
\textbf{ii:} $\forall y\in$ found$(T)$,  $ \
\lambda_{xy}=2^{t-1}-2^{t-u-1}$.\\
 In this case the non-existence of T can be proved the same way as in
the case $t=i+2$ (Case A.2(b)ii). }

\end{proof}
\vspace{.2cm} \textbf{Remark:} In this sequel we improved the
conjecture for Steiner trade. It remained for interested reader
to prove the correctness of this conjecture for any trade.
\vspace{.2cm}\\
 \textbf{Acknowledgement:} The authors thank the
referees who helped us to improve  presentation of this
manuscript.

\end{document}